\documentclass[amsfonts]{article}
\usepackage{graphicx}

\newtheorem{proposition}{Proposition}
\newtheorem{definition}[proposition]{Definition}
\newtheorem{lemma}[proposition]{Lemma}
\newtheorem{remark}{Remark}
\newtheorem{theorem}{Theorem}
\newcommand{\qedsymb}{\mbox{ }~\hfill~{\rule{2mm}{2mm}}}
\newenvironment{pf}{\begin{trivlist}
\item[\hspace{\labelsep}{\bf \noindent Proof. }]
}{\qedsymb\end{trivlist}}

\begin{document}

\title{Orthomodular Lattices and Quantales}

\author{Leopoldo Rom\'{a}n}

\maketitle

\begin{center}
{\it Dedicated to Raquel Hern\'{a}ndez}
\end{center}

\bigskip
\noindent {\bf Abstract.}
Let $L$ be a complete orthomodular lattice. There is a one to one
correspondence between complete boolean subalgebras of $L$ contained
in the center of $L$ and endomorphisms $j$ of $L$ satisfying the
Borceux-Van den Bossche conditions.

\bigskip

\section*{Introduction}

In [9], we study the notion of an idempotent right-sided quantale
versus the concept of orthomodular lattice. We proved that a complete
orthomodular lattice $L$ has a natural structure of an idempotent
right-sided quantale if we take the central cover of an arbitrary
element $a$ (denoted by $e(a) )$ of $L$; this construction induces an
endomorhism of $L$ satisfying certain conditions, see section $I$ for
more details. Therefore, there is a natural question after this
claim: which endomorphisms of $L$ can produce an idempotent,
right-sided quantale?

Recall that when the endomorphism $ j:L \rightarrow L$ satisfies
the Borceux-Van den Bossche's conditions (denoted by B.-V.B.), then
$j$ induces an idempotent, right-sided quantale in $L$; in fact, the
structure is quite simple: if $a,b$ are arbitrary elements of $Q$
then $a\&b = a \wedge j(b)$.  The reader can see [3] or [9] for
details. The endomorphism $j$ is a closure operator satisfying three
conditions which are related with the concepts of nucleus and quantic
nucleus considered in intuitionistic logic and quantum logic. In [8],
Beatriz Rumbos and myself gave a characterization of nuclei in
orthomodular lattices and quantic lattices which in particular
produces a characterization of quantic nuclei in orthomodular
lattices. The difference here is the binary operation $\&^F$
considered in [8]. This operation was first considered by P.D. Finch
in [4], where he suggested that $\&^F$ has very similar properties
with the connectives $( \wedge, \rightarrow )$, considered in
classical logic and intuitionistic logic.In fact, this is true but
the main difference is the lack of the associativity property of
$\&^F$. Nevertheless, $-\&^F b : L \rightarrow L$ has a right adjoint
which is the Sasaki hook, for an arbitrary element $b$ of $L.$
Unfortunately, the operation $b \&^F - :L \rightarrow L$ does not
have a nice property such as: $\&^F$ is associative, $b \&^F -$ has a
right adjoint or $b \&^F -$ preserves order. In fact, if $L$
satisfies one of these conditions then $L$ will be a boolean algebra
and conversely, if $L$ is a boolean algebra then $a\&^F b$ is just $a
\wedge b$. The interested reader can consult [8] for more details.

After all these comments, I would like to say some words about the
result I will prove. Last year, Prof. R. Greechie suggested the
following idea: perhaps the work of M.F. Janowitz can help to
characterize all the endomorphisms of a complete orthomodular lattice
satisfying the B.-V.B.  conditions; I must say he was quite right,
the paper due to M.F. Janowitz, under the title ``Residuated Closure
Operators'' gave me some ideas which are crucial in the proof of the
main result of this article. See [6] for more details.

 Finally, the reader must note that by ``the logic of quantum mechanics'' we
mean the lattice theoretic ``quantum logics'' of Birkhoff and von
Neumann [2], hence we do not consider the quantales introduced
by David N. Yetter in [10], under the name of ``Girard Quantales'' where
he considers a different logic for quantum mechanics; roughly speaking
the logic considered by David Yetter is a logic involving an asociative
(in general noncommutative) operation  ``and then''. Yet, Girard quantales
have a close relation with linear logic, a logic introduced by
J.Y. Girard, the reader can consult [4] for details and comments
about this logic. Clearly, a natural question is if linear logic
has some relation with the quantum logic introduced by Birkhoff
and von Neumann, we will look about this problem in a future work.

The article is organized as follows. In the first section we
introduce the concepts we need for our purposes. In the second
section we show the main results of this article. We prove that if we
have an arbitrary $ j:L \rightarrow L$ satisfying the Borceux-Van den
Bossche conditions then $j(L)$ is a boolean subalgebra of $L$
contained in the center of $L$ and conversely every boolean
subalgebra $M$ of $L$ contained in the center of $L$ induces an
endomorphism $k$ satisfying the B.-V.B. conditions. Hence there is a
one to one correspondence between endomorphisms $j$ satisfying the
B-V.B. conditions and boolean subalgebras of $L$ contained in the
center of $L$.

I want to express my sincere thanks to Prof. R. Greechie for his
suggestion.

\begin{center}
\textbf{Section I}
\end{center}

\begin{definition}
A quantale $Q$ is a lattice having arbitrary joins $\vee $ together
with an associative product $\&$ such that:

\begin{enumerate}
\item  $a\&(\vee _{i\in I}b_i)=\vee _{i\in I}(a\&b_i);$

\item  $(\vee _{i\in I}a_i)\&b=\vee _{i\in I}(a_i\&b)$
\end{enumerate}

for all $a,b,a_i,b_i\in Q.$

Moreover, we will say that the quantale $Q$ is an idempotent and right-sided
if it satisfies the following two conditions

\begin{enumerate}
\item  $a\&1=a;$

\item  $a\&a=a$ for all $a\in Q.$
\end{enumerate}
\end{definition}

\begin{remark}
In $\left[ 3\right] $ F.Borceux and G. Van Den Bossche proved that
given any complete lattice $(Q,\leq )$ there is a one-to-one
correspondence between binary operations $\&:Q\times Q\rightarrow Q$
making $Q$  an idempotent, right-sided quantale and closure
operations $j:$ $Q\rightarrow Q$ satisfying the following axioms
(these are the B.-V.B. conditions):

\begin{enumerate}
\item  $a\leq j(a);$

\item  $j(a\wedge j(b))=j(a)\wedge j(b);$

\item  $a\wedge j(\vee _{i\in I}b_i)=\vee _{i\in I}(a\wedge j(b_i));$

\item  $(\vee _{i\in I}a_i)\wedge j(b)=\vee _{i\in I}(a_i\wedge j(b)).$
\end{enumerate}

\end{remark}

 We just mention a simple consequence of this result. If we have an
endomorphism $j$ of a complete lattice $Q$ satisfying the B.-V.B.
conditions and we take the fixed points of $j$ (denoted by $Q_j$) then
$Q_j$ is not only a quantale is a locale as the reader can ckeck easily.
Hence, in logical terms any idempotent and right-sided quantale $Q$
has a locale as a sublattice;i.e., $Q$ has a model of intuitionistic
logic.

In the rest of this article we shall work only with an idempotent
right-sided quantale. There are many examples of these quantales. For
instance, any Locale $H$ is an idempotent, right-sided quantale in a
trivial way; the binary operation $\&$ is just $\wedge$. The closed
ideals of a $C^{*}$-algebra is also an example, here the binary
operation $\&$ is just the closure of the product of two ideals. For
the next example we need some definitions. We begin with the following

\begin{definition}
A complete lattice $L=(L,\vee ,\wedge ,\perp )$ is a complete
orthomodular lattice if there exists a unary operation
$\bot:L\rightarrow L$ satisfying the conditions:

\begin{enumerate}
\item $a^{\bot\bot} =a.$

\item $(\bigvee_{i \in I}a_{i})^\bot = \bigwedge_{i \in I}a_{i}^\bot,
$ for all $a,a_{i} \in L$ and any set $I$

\item $ a \vee a^\bot = 1.$

\item $ a \wedge a^\bot = 0.$
\end{enumerate}

Moreover, $L$ satisfies  the following weak modularity property:

\noindent Given any $a,b$ $\in L$ with $a\leq $ $b$ then $b=a\vee
(a^{\perp }\wedge b)$ (equivalently $a=(a\vee b^{\perp })\wedge b$).
\end{definition}

As we said in the introduction, if we have a complete orthomodular
lattice $L$, the way of inducing an idempotent, right-sided quantale
structure in $L$ is by taking the {\it{central cover of an element}}.
We shall introduce more concepts.

First of all, given two arbitrary elements $a,b$ of an arbitrary
complete orthomodular lattice $L$, $a\&^F b = ( a \vee b^\bot )
\wedge b$ and the Sasaky hook is given by the following rule:
$a\rightarrow b=(a\wedge b)\vee a^{\bot }$. Since we are interesed in
the center of a complete orthomodular lattice, we consider first the
notion of compatibility:

\begin{definition}
Let $L$ be a complete orthomodular lattice. We say $a,b\in L$ are
compatible elements (denoted by $bCa$ ) if and only if
$b\&^Fa=a\wedge b.$
\end{definition}

Notice that whenever $a,b$ are compatible elements it is easy to see
that $a\&^Fb=a\wedge b$ also holds.

The simplest example of a pair of elements $a,b$ which are compatible
is whenever one of these elements belongs to the center $Z(L)$ of the
complete orthomodular lattice $L$. The definition of the center is as
follows:

\begin{definition}
Let $L$ be an arbitrary complete orthomodular lattice. The center of
$L$, denoted by $Z_F(L)$) is the set
$$Z_F(L)=\left\{ a\in L\mid a\&^Fb=b\&^Fa=a\wedge b\mbox{\ for all\ }b\in
L\right\} $$
\end{definition}

Notice that $Z_F(L)$ is a boolean subalgebra of $L$. In particular,
$0,1$ belong always to the center of $L$. Hence $Z_F(L)$ is
non-empty. We define now the central cover of an arbitrary element of
$L$.

\begin{definition}
Let $L$ be a complete orthomodular lattice. If $a$ is an arbitrary
element of $L$, the central cover of $a$ is given by:

\begin{center}
$e(a)=\wedge \left\{ z\in Z_F(L)\mid a\leq z\right\} .$
\end{center}
\end{definition}

The central cover of an element always exists since $1$ is an element
of this set. The reader can see [1], p.129 and also the comments
contained in that book. The next proposition give us the example we
are interested in.

\begin{proposition}$[9]$
The map $e:L\rightarrow L$ satisfies the Borceux-Van Den Bossche
conditions. $L$ has a binary operation $\&$, defined by $a\&b = a
\wedge e(b)$, making it an idempotent, right-sided quantale.
\end{proposition}

Therefore, this is the first example of an endomorphism $j$ of $L$
making it an idempotent, right-sided quantale. Clearly, the
orthomodular lattice must be complete if one wants to preserve the
definition of a quantale. However, almost all the concepts described
above, can be defined in an arbitrary orthomodular lattice. If we
take the classical example of the closed subspaces of a Hilbert space
$H$, the center is trivial; the only elements of the center are
$0,1$. In the literature an orthomodular lattice is called
irreducible whenever it has trivial center. Clearly, if we have an
irreducible orthomodular lattice the quantale structure that we get
is not really interesting. However, this phenomenon does not occur
always.  We just mention one example of a finite orthomodular lattice
with non-trivial center. Namely, $G_{12}$.

\begin{center}
\includegraphics{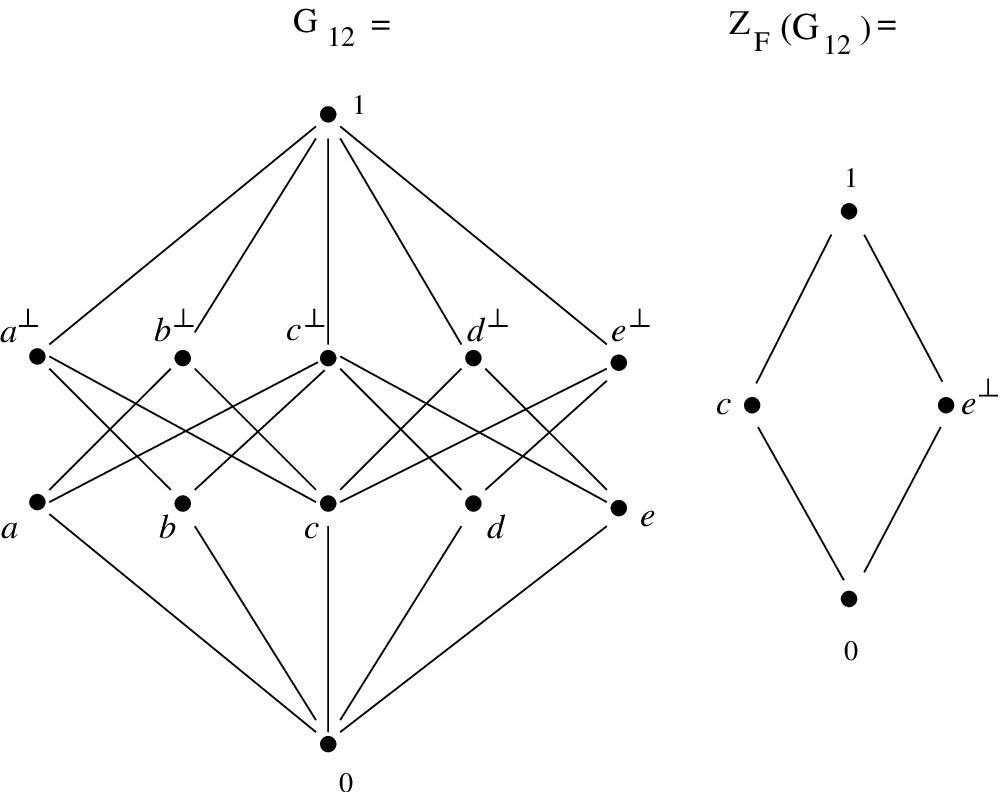}
\end{center}

We close this section with another comment. The binary operation
$\&^F$ is really important for the construction of the second binary
operation $\&$; in fact, all the concepts we had were defined in
terms of $\&^F$, despite the lack of associativity or equivalently
that the endomorphism $a \&^F - $ does not necessarily preserves
order.

\begin{center}
\textbf{Section II}
\end{center}

We shall prove now the main result of this article. We will assume
$L$ is an arbitrary complete orthomodular lattice. If we take an
arbitrary endomorphism $j$ of L satisfying the B.-B.V. conditions
then clearly $L$ is an idempotent, right-sided quantale. We just
define $ a \& b = a \wedge j(b) $, for elements $a,b$ in $L$. Now,
the question is: which endomorphisms $j$ of $L$ satisfy the B.-V.B.
conditions? Actually, can we give a characterization of such
endomorphisms in terms of another concept? We shall see that this is
the case. We would like to recall that some of the results are
inspired by the work of M.F. Janowitz contained in $[6]$. First of
all, we begin with the following

\begin{lemma}
If $ j:L \rightarrow L$ is an arbitrary endomorphism of a complete
orthomodular lattice, then given any element $a$ of $L$, $j(a)$
satisfies the following identity.

\begin{center}
$j(a)=\wedge \left\{ x\in L\mid j(x)=x, a\leq x\right\} .$
\end{center}
\end{lemma}

\begin{pf}
Indeed, let us call $z$ the RHS of the last equality. Since, $a \leq
j(a)$ and $j$ is idempotent, we have $z \leq j(a)$. Now, $a \leq z$
and since $j$ preserves order we get: $ j(a) \leq j(z) = z$. Hence,
$z = j(a)$. As we claimed.
\end{pf}

We must notice that $j(0)$ is equal to $0$. The reason is quite
simple, just take the empty set for $I$ in the third property of
B.-V.B. conditions. We now prove the next

\begin{proposition}
Suppose $L$ is an arbitrary complete orthomodular lattice and $j$
is an endomorphism satisfying the B.-V.B. conditions then the subset
$L_{j}$ of $L$ defined by
$$L_{j}= \left\{ x\in L\mid j(x) = x \right\} .$$
is a complete boolean sublattice of $L$ contained in the center of
$L$.
\end{proposition}

\begin{pf}
We check first, $L_{j}$ is a complete lattice. If $\{b_i\}_{i \in I}$
is an arbitrary family of elements of $L_{j}$ and taking $a=1$, by
the third property of the B.-V.B. conditions we have:
$$j(\vee _{i\in I}b_i) = 1 \wedge j(\vee _{i\in I}b_i)=\vee _{i\in
I}(1 \wedge j(b_i))=\vee _{i\in I}j(b_i)=\vee _{i\in I}b_i.$$

Hence, $L_{j}$ is a complete lattice. We shall see now $L_{j}$ is
closed under complements; i.e., if $a \in L_{j}$ then $a^\bot$ also
is an element of $L_{j}$.  We already knew $a^\bot \leq j(a^\bot)$.
We only need to check: $j(a^\bot) \leq a^\bot$. This is equivalent
to: $j(a^\bot) \&^F a \leq 0$. Since $ - \&^F a$ has a right adjoint,
namely the Sasaky hook, it is enough to check this. We calculate
$j(a^\bot) \&^F a$.
$$j(a^\bot) \&^F a = (j(a^\bot) \vee a^\bot) \wedge a = a \wedge j(a^\bot).$$

Since $ a^\bot \leq j(a^\bot)$. Now, since $a \in L_{j}$, we shall
see, $a \wedge j(a^\bot) = 0$.
$$a \wedge j(a^\bot) =j(a )\wedge j(a^\bot) = j(a^\bot \wedge j(a))
\leq j(a^\bot \&^F j(a)) = j(a^\bot \& a)= j(0)=0$$

Hence, $a \wedge j(a^\bot) =0$ and $j(a^\bot) \leq a^\bot$. Therefore
$L_{j}$ is closed under complements.

We check now $L_{j}$ is a boolean sublattice of $L$. Indeed, if
$a,\{b_i\}_{i \in I}$ belong to $L_{j}$, from the third property of
the B-V.B. conditions we have:
$$a \wedge (\vee _{i\in I}b_i) = a \wedge j(\vee _{i\in I}b_i) =\vee
_{i\in I}(a \wedge j(b_i)) =\vee _{i\in I}(a \wedge b_i).$$

Therefore $L_{j}$ is a distributive lattice closed under complements;
i.e., $L_{j}$ is a boolean algebra. Finally, we shall see $L_{j}$ is
contained in the center of $L$.

Suppose $a$ and $j(b)$ are arbitrary elements of $L$ and $L_{j}$
respectively.  We calculate $a \&^F j(b)$.
$$a \&^F j(b) = ( a \vee j(b)^\bot ) \wedge j(b) = [a \vee j(b^\bot)]
\wedge j(b) = ( a \wedge j(b) ) \vee ( j(b)^\bot \wedge j(b) ) = a
\wedge j(b).$$

Since $L_{j}$ is closed under complements and by the third property
of the B.-V.B. conditions. In a similar way, we can check $ j(b) \&^F
a = j(b) \wedge a$. Hence $j(b)$ is in the center of $L$.
\end{pf}

\begin{remark}

  In the proof of the proposition we used the inequality $ a^\bot \& a \leq
a^\bot \&^F j(a)$. Actually, it is easy to see that $a \& b \leq a \&^F j(b)$
for arbitrary elements $a,b$ in $L$. Moreover, $j(L)$ is not only a boolean
subalgebra of the center of $L$, it is a complete boolean subalgebra of the
center of $L$.
\end{remark}

 We will see now the converse of this proposition.

\begin{proposition}

 Let $L$ be an arbitrary complete orthomodular lattice. If $Z_F(L)$ denotes the
center of $L$ and $M$ is a complete subalgebra of $L$ then the endomorphism
$j:L \rightarrow L$ defined by

\begin{center}
$j_M(a)=\wedge \left\{ x\in M \mid a\leq x\right\} .$
\end{center}

\noindent satisfies the B-V.B. conditions and therefore $L$ is an idempotent, right-sided quantale. The binary operation $\&$ is given by $ a \& b = a \wedge j_M(b) $ where $a,b$ are arbitrary elements of $L$.
\end{proposition}

\begin{pf}
 Clearly, given any element $a$ of $L$ we have: $a \leq j_M(a)$ and $j_M$ is
idempotent. Moreover, if $z \in M$ then $j_M(z) = z$; using these results, $j_M$
preserves order as the reader can check easily. Now, we shall see $j_M$ preserves arbitrary suprema. Suppose $\{a_i\}_{i \in I}$ is a family of elements of $L$. Since $j_M$ preserves order, we only need to verify
$j_M(\vee_{i \in I}a_i)
\leq \vee_{i\in I}j_M(a_i)$ but this is true since $M$ is a complete subalgebra of $Z_F(L)$ and therefore $ \vee_{i\in I}j_M(a_i) \in M$.

 It is not hard to prove the third and the fourth properties of the B.-V.B. since $j_M(a)$ is an element of the center of $L$ and $j_M$ perserves arbitrary  joins. We only need to verify the second property; i.e., $j_M(a \wedge j_M(b))=
j_M(a) \wedge j_M(b)$. Given any element $z$ of $M$ and an arbirary element $a$
of $L$ the identities hold:

\begin{center}
$ z \wedge j_M(a) = \{ z \wedge j_M(z \wedge a) \} \vee \{ z \wedge j_M(z^\bot
\wedge a) \} = \{ z \wedge j_M(z \wedge a) \} = j_M(z \wedge a)$
\end{center}

 In particular, $j_M(a \wedge j_M(b))= j_M(a) \wedge j_M(b)$. Hence, $j_M$ satisifes the B.-V.B. conditions as we claimed.

\end{pf}

We summarize the results in the following

\begin{theorem}

        Let $L$ be a complete orthomodular lattice. There is a one to one
correspondence between complete boolean subalgebras $M$ contained in the center
of $L$ and endomorphisms $j:L \rightarrow L$ satisfying the Borceux-Van den
Bossche conditions.
\end{theorem}

        We would like to mention another example of an orthomodular lattice
with trivial center. The lattice is called $M_{On}$ for $ 1 \leq n $. In fact,
$M_{On}$ is not only orthomodular, it is modular but it is not distributive.

\begin{center}
\includegraphics{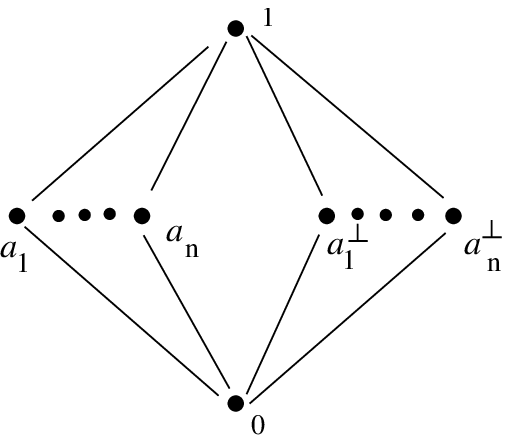}
\end{center}

         It is not hard to check that given $ i \neq j $,  $a_i \&^F a_j = a_j$ and  $a_j \&^F a_i = a_i$. Hence, the center is trivial.

E-mail address: leopoldo@matem.unam.mx.\\[5pt]
Instituto de Matem\'{a}ticas, UNAM\\
Area de la Investigaci\'{o}n Cient\'{\i}fica, Ciudad Universitaria\\
04510, M\'{e}xico, D.F.\\
Mexico


\begin{thebibliography}{99}
\bibitem{1} Beltrametti Enrico G. and Cassinelli Gianni, The logic of
Quantum Mechanics, Encyclopedia of Mathematics and its applications, $15$, $
1981$.

\bibitem{2} Birkhkoff G. and von Neumann J., The logic of quantum mechanics,
Annals of Mathematics, ser.2, vol.$37$,823-246, $1936$.

\bibitem{3} Borceux F. and Van Den Bossche G., Quantales and their
sheaves, Order $3$, $1986$, 61--87.

\bibitem{4} Finch P. D, Quantum logic as an implication algebra,
Bull. Aust. Math. Soc., 1, 101--106, $1970$.
\bibitem{5} Girard J.Y., Linear Logic, Theoretical Computer Science,
50, 1-102,1987.

\bibitem{6} Janowitz M.F., Residuated Closure Operators, Portugal Math., 26,
1967, 221-252.

\bibitem{7} Rom\'{a}n L. and Rumbos B., A characterization of nuclei in
orthomodular lattices and quantic nuclei, Journal of Pure and Applied
Algebra, $73$, 155-163, $1991.$

\bibitem{8} Rom\'{a}n L. and Rumbos B., Quantum Logic Revisted,
Foundations of Physics, $21$ No. 6, 727--734, $1991.$

\bibitem{9} Rom\'{a}n Leopoldo and Zuazua Rita, Gelfand Quantales and
Orthomodular Lattices. Submmitted to Journal of Pure and Applied Algebra.


\bibitem{10} Yetter N.D., Quantales and (Non Commutative) Linear Logic,
The Journal of Symbolic Logic, $55$, No.1, 41-64,1990.

\end{thebibliography}
\end{document}